\documentclass{birkjour}
\usepackage{amsmath,epsfig,amssymb,graphics}
\usepackage{caption}
 \newtheorem{thm}{Theorem}[section]

 \theoremstyle{definition}
 
 \theoremstyle{remark}
 \newtheorem{rem}[thm]{Remark}
 
 \numberwithin{equation}{section}

\begin{document}

%
%
%
%
%
%
%
%
%

\title[A new analytic method for solving nonlinear problems]
 {A new analytic method with a convergence-control parameter for solving nonlinear problems}

\author[X.L. Zhang]{Xiaolong Zhang}

\address{%
School of Mathematical Sciences\\
Dalian University of Technology\\
116024\\
China}

\email{topmaths@mail.dlut.edu.cn}

\thanks{This work was partially supported by the National Natural Science Foundation of China (51379033, 51522902).}
\author[S.X. Liang]{Songxin Liang}
\address{School of Mathematical Sciences\\
Dalian University of Technology\\
116024\\
China}
\email{sliang668@dlut.edu.cn}


\date{September 1, 2016}

\begin{abstract}
In this paper, a new analytic method with a convergence-control parameter $c$ is first proposed.
The parameter $c$ is used to adjust and control the convergence region and rate of the resulting series solution.
It turns out that the convergence region and rate can be greatly enlarged by choosing a proper value of $c$.
Furthermore, a numerical approach for finding the optimal value of the convergence-control parameter is given.
At the same time, it is found that the traditional Adomian decomposition method is only a special case of the new method.
The effectiveness and applicability of the new technique are demonstrated by several physical models including
nonlinear heat transfer problems, nano-electromechanical systems, diffusion and dissipation phenomena, and dispersive waves.
Moreover, the ideas proposed in this paper may offer us possibilities to greatly improve current analytic and numerical techniques.
\end{abstract}

\maketitle

\section{Introduction}\label{sect:intro}
Nonlinear problems occur in almost every field in science and engineering. The exact solutions are very difficult to obtain for most strong nonlinear differential equations. To solve these problems, several analytic methods such as the traditional perturbation method \cite{Nayfeh2000}, Lyapunov's artificial small parameter method \cite{Lyapunov1992}, the $\delta-$expansion method \cite{Karmishin1990} and the Adomian decomposition method (ADM) \cite{Adomian1983} have been proposed.

The Adomian decomposition method, which was first proposed by George Adomian in the 1980s, is a powerful approach for seeking series solutions to linear and nonlinear problems. With the development of computer softwares such as Maple and Mathematica, the advantages of the method for solving complicated nonlinear problems are more and more obvious. It has been successfully applied to solve many problems in physical sciences such as the nonlinear Klein-Gordon equation \cite{Deeba1995}, the Lane-Emden problem \cite{Wazwaz2014}, the hydromagnetic peristaltic flow \cite{Ebaid2008},
the nonlinear fin problem \cite{Singla2015} and the plate flow \cite{Siddiqui2010}. To make the method more effective, many researchers have done lots of improvements, for example, giving the noise terms \cite{Adomian1992}, comparing with other methods \cite{Rach1987,Kaya2004}, discussing the convergence rate \cite{Wazwaz1999}, analysing the error \cite{El-Kalla2011}, modifying the recursive scheme \cite{Duan2011,Lin2014}and indicating the essence \cite{Zhang2015}. Recently, Duan and his colleagues introduced a convergence parameter into the ADM firstly \cite{Duan2010,Lu2014,Duan2013}.

For many real-world problems, one expects larger convergence regions and rates for the resulting series solutions. However,
the traditional ADM usually cannot reach the goal, as shown in Section \ref{subsect:4.1} below. In the present paper, we introduce a new more general method called Adomian Decomposition Method with a convergence-control Parameter (ADMP), which can be used successfully to achieve the goal and is different with various modifications of the ADM by other researchers. Furthermore, the method offers a family of solutions depending on the convergence-control parameter $c$,
providing a chance to choose the optimal series solution. It turns out that the convergence region and rate of the resulting series solution
can be greatly enlarged by choosing a proper value of $c$.

In Section \ref{sect:TADM}, we present a brief review of the traditional ADM. In Section \ref{sect:ADMP}, the ADMP is described explicitly.
Moreover, a numerical approach for finding the optimal value of the convergence-control parameter is also given. To demonstrate the effectiveness
and applicability, several nontrivial applications including nonlinear ODEs and nonlinear PDEs are discussed in Section \ref{sect:applica}.
\section{Basic ideas of the traditional ADM}\label{sect:TADM}

We shall begin to consider the nonlinear differential equation
\begin{equation}\label{equ:pro}
L[u]+R[u]+N[u]=f(t),\;t\in\Omega
\end{equation}
with some initial/boundary conditions, where $u=u(t), L, R, N$ and $f(t)$ denote some unknown function, an easily invertible linear operator
(usually the highest-order linear operator), the remainder of the linear operator, the nonlinear operator and the given source term.
The linear operator $L$ is designed for performing the inverse easily, but the choice is not unique. For example, for the Lane-Emden equation
$u''(t)+\frac{p}{t}u'(t)+f(u(t))=0$, one can choose either the highest-order differential operator $\frac{d^{2}}{dt^{2}}(\cdot)$ or $t^{-p}\frac{d}{dt}(t^{p}\frac{d}{dt}(\cdot))$
as $L$. The nonlinear operator $N$ is analytic. Multiplying (\ref{equ:pro}) by $L^{-1}$, one obtains
\begin{equation}\label{equ:inv}
L^{-1}[L(u)]=L^{-1}[f(t)]-L^{-1}[R(u)]-L^{-1}[N(u)].
\end{equation}
According to the definition of integral operator, one has
\begin{equation}\label{equ:inte}
L^{-1}[L(u)]=u-\phi(t),
\end{equation}
where $\phi(t)$ is determined by the initial/boundary conditions mentioned before.
Combining Eq. (\ref{equ:inv}) with Eq. (\ref{equ:inte}), one gets
\begin{equation}
u=L^{-1}[f(t)]+\phi(t)-L^{-1}[R(u)]-L^{-1}[N(u)].
\end{equation}

In the traditional ADM, the solution $u(t)$ is expressed as the decomposition series
\begin{equation}\label{equ:seri}
u(t)=\sum_{k=0}^{+\infty}u_{k}(t),
\end{equation}
where $u_k(t),k\ge 0$ are determined by the recursive scheme
\begin{equation}\label{equ:part}
\begin{aligned}
&u_0(t)=L^{-1}[f(t)]+\phi(t),\\
&u_k(t)=-L^{-1}[R(u_{k-1}(t))]-L^{-1}[A_{k-1}],\; k\ge1,
\end{aligned}
\end{equation}
and
\begin{equation}\label{equ:poly}
A_k=\frac{1}{k!}\left.\left[\frac{d^{k}}{d\epsilon^{k}}N\left(\sum_{i=0}^{+\infty}u_i(t)\epsilon^{i}\right)\right]\right|_{\epsilon=0}, \; k\ge0
\end{equation}
are decomposition polynomials.

To specify the decomposition polynomials in Eq. (\ref{equ:poly}), one has
\begin{equation*}
\begin{aligned}
A_0&=N[u_0(t)],\\
A_1&=N'[u_0(t)]u_1(t),\\
A_2&=N'[u_0(t)]u_2(t)+\frac{1}{2!}N''[u_0(t)]u_1^{2}(t),\\
A_3&=N'[u_0(t)]u_3(t)+N''[u_0(t)]u_1(t)u_2(t)+\frac{1}{3!}N'''[u_0(t)]u_1^{3}(t),\\
&\vdots
\end{aligned}
\end{equation*}
It is seen that $A_k$ depends on $A_i, i=0, 1, \cdots, k-1$ and the nonlinear operator $N$.
In practice, instead of finding a solution expression (\ref{equ:seri}), one can only calculate an $n$th-order approximation
\begin{equation}\label{equ:appro}
\psi_n(t)=\sum_{k=0}^{n}u_k(t),
\end{equation}
of which $u_k(t), k=0, 1, \cdots, n$ can be calculated recursively via (\ref{equ:part}).

\section{Description of the new ADMP}\label{sect:ADMP}

Based on the traditional ADM, a convergence-control parameter $c$ and a artificial parameter $\epsilon$ are introduced to Eq. (\ref{equ:pro}),
so a new equation is
\begin{equation}\label{probp}
L[u]+\left(\epsilon c+\epsilon^2(1-c)\right)\left(R[u]+N[u]\right)=f(t).
\end{equation}
One then set
\begin{equation}\label{serip-epsi}
u(t)=\sum_{k=0}^{+\infty}v_{k}(t,c)\epsilon^k.
\end{equation}
Multiplying (\ref{probp}) by $L^{-1}$, one has
\begin{equation}\label{equ:prob-inv}
u(t)=-\left(\epsilon c+\epsilon^2(1-c)\right)L^{-1}\left(R[u]+N[u]\right)+L^{-1}[f(t)]+\phi(t).
\end{equation}
Substituting (\ref{serip-epsi}) into (\ref{equ:prob-inv}) and collecting the coefficients of powers of $\epsilon$, a new recursive scheme is obtained
\begin{equation}\label{equ:partp}
\begin{aligned}
v_0(t,c)    =&L^{-1}\left[f(t)\right]+\phi(t),\\
v_1(t,c)    =&-cL^{-1}\left[R(v_0(t,c))+B_0\right],\\
v_{k}(t,c)=&-cL^{-1}\left[R(v_{k-1}(t,c))+B_{k-1}\right]\\
&-(1-c)L^{-1}\left[R(v_{k-2}(t,c))+B_{k-2}\right],k\ge2
\end{aligned}
\end{equation}
where the Adomian polynomials are given by
\begin{equation}\label{equ:polyp}
B_k=\frac{1}{k!}\left.\left[\frac{d^{k}}{d\epsilon^{k}}N\left(\sum_{i=0}^{+\infty}v_i(t,c)\epsilon^{i}\right)\right]\right|_{\epsilon=0}, \; k\ge0.
\end{equation}
To specify the polynomials in Eq. (\ref{equ:polyp}), one has
\begin{equation*}
\begin{aligned}
B_0&=N[v_0(t,c)]=N[u_0(t)]=A_0,\\
B_1&=N'[v_0(t,c)]v_1(t,c)=cN'[u_0(t)]u_1(t)=cA_1,\\
B_2&=N'[v_0(t,c)]v_2(t,c)+\frac{1}{2!}N''[v_0(t,c)]v_1^{2}(t,c)=c^{2}A_2+(1-c)A_1,\\
B_3&=N'[v_0(t,c)]v_3(t,c)+N''[v_0(t,c)]v_1(t,c)v_2(t,c)+\frac{1}{3!}N'''[v_0(t,c)]v_1^{3}(t,c)\\
&=c^{3}A_3+2c(1-c)A_2,\\
&\vdots
\end{aligned}
\end{equation*}

The relationship between the new recursive scheme (\ref{equ:partp}) and the traditional recursive scheme (\ref{equ:part}) is as follows.
\begin{equation*}
\begin{aligned}
v_0(t,c)&=u_0(t),\\
v_1(t,c)&=cu_1(t),\\
v_2(t,c)&=c^{2}u_2(t)+(1-c)u_1(t),\\
v_3(t,c)&=c^{3}u_3(t)+2c(1-c)u_2(t),\\
\vdots
\end{aligned}
\end{equation*}
When $\epsilon=1$ in (\ref{serip-epsi}), the solution $u(t)$ with convergence-control parameter $c$ is presented by
\begin{equation}\label{equ:serip}
v(t,c)=\sum_{k=0}^{+\infty}v_k(t,c).
\end{equation}
Similar to Eq. (\ref{equ:appro}), one can only calculate an $n$th-order approximation
\begin{equation}\label{equ:approp}
{\psi}_n(t,c)=\sum_{k=0}^{n}v_k(t,c).
\end{equation}

Next, an approach for determining the optimal value of the convergence-control parameter $c$ for the $n$th-order approximation is given.
The optimal $c$ is determined by solving the equation
\begin{equation}\label{equ:opti}
\frac{\partial\mathbf{E}(c)}{\partial c}=0,
\end{equation}
where the averaged residual error $\mathbf{E}(c)$ of the approximation (\ref{equ:approp}) on $\Omega$ is defined by
\begin{equation*}
\mathbf{E}(c)=\frac{1}{M}\sum_{k=1}^{M}\left[\left(L+R+N\right)({\psi}_n(t_{k},c))-f(t_k)\right]^{2},
\end{equation*}
and $t_1, t_2, \cdots,t_M \in \Omega$ are sample points.

\begin{rem}
It is easily seen that the traditional ADM is only a special case of the new approach ADMP when $c=1$.
\end{rem}

\section{Applications}\label{sect:applica}

\subsection{A nonlinear heat transfer problem} \label{subsect:4.1}

Consider the nonlinear heat transfer problem governed by a nonlinear ordinary differential equation \cite{Liang2009}
\begin{equation}\label{appl1:prob}
\left(1+\epsilon u(t)\right)u'(t)+u(t)=0,\; u(0)=1,
\end{equation}
where $\epsilon$ is a physical parameter, and the prime denotes differentiation with respect to the time $t$. The closed-form solution is unknown, but one can get
\begin{equation}\label{appl1:der-exact}
 u'(0)=-\frac{1}{1+\epsilon}.
\end{equation}

In the following, the traditional perturbation method, the traditional ADM and the new ADMP
are applied in sequence to seek series solutions to $u'(0)$, and comparisons between them are made.

\begin{figure}[ht]
\setcaptionwidth{2in}
\centering
\begin{tabular}{cc}
\begin{minipage}{2in}
\includegraphics[width=2in,height=2in]{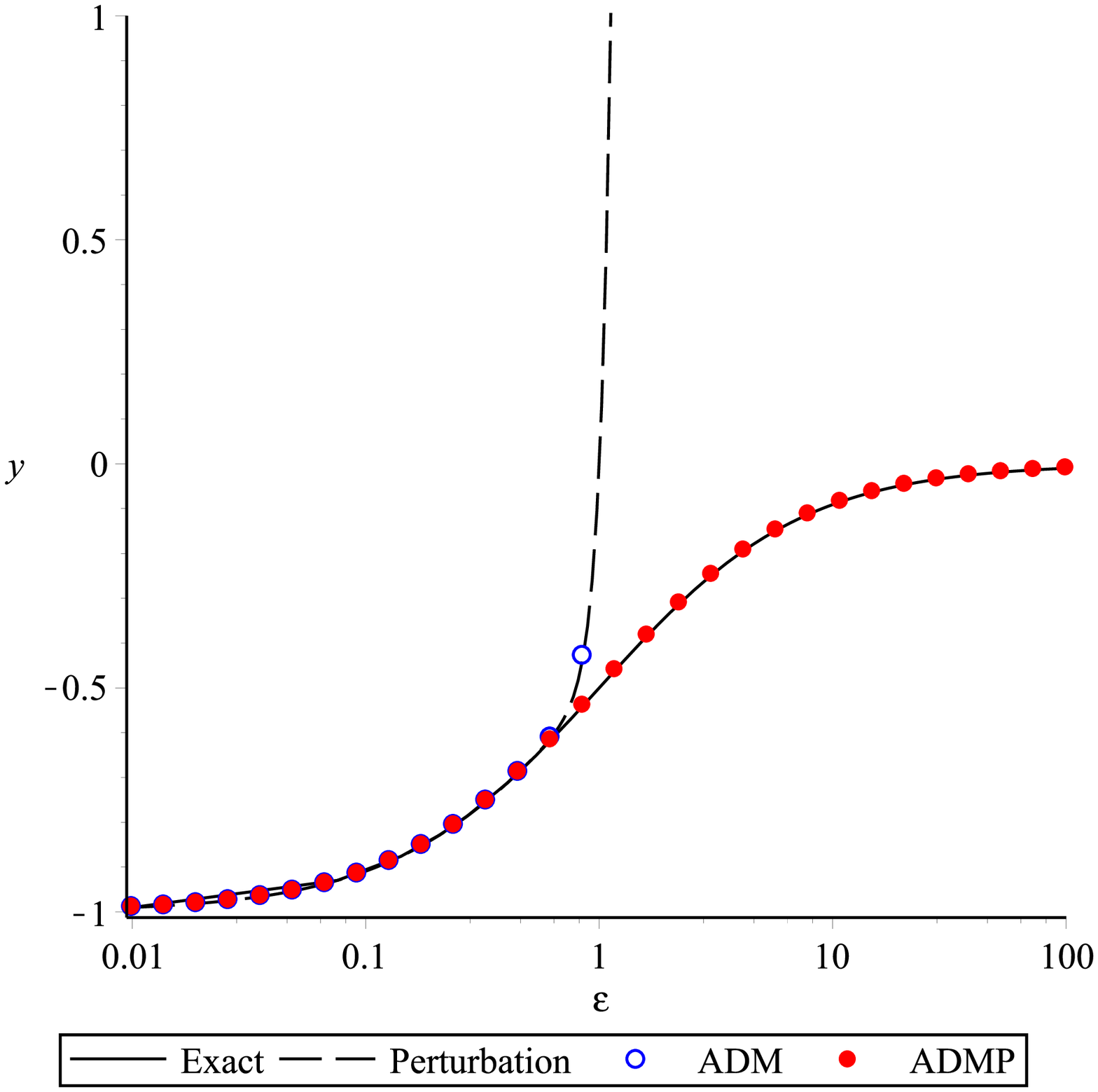}
\caption{Solid line: exact solution (\ref{appl1:der-exact}); dash line: $9$th-order perturbation approximation; circle line: $10$th-order ADM; solid circle line: $1st$-order ADMP with $c=\frac{1}{1+\epsilon}$.}\label{fig:appl1-1}
\end{minipage}
\hspace{0.1in}
\begin{minipage}{2in}
\includegraphics[width=2in,height=2in]{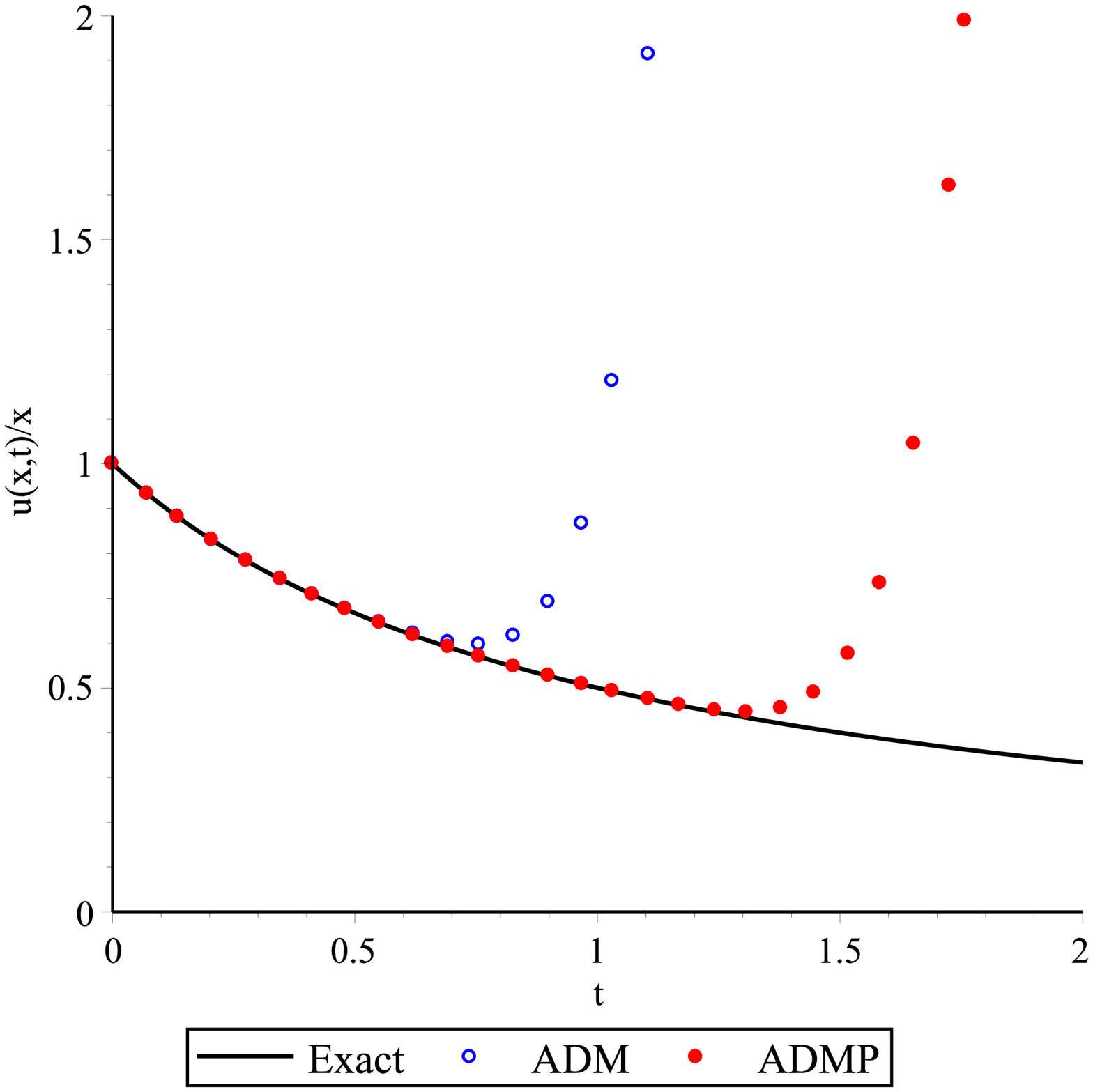}
\caption{ Solutions for (\ref{appl3:problem}). Solid line: $\frac{u(x,t)_{ex}}{x}$; circle line:$\frac{ \psi_{10}(x,t,1)}{x}$; solid circle line : $\frac{\psi_{10}(x,t,0.8)}{x}$.}\label{fig:appl3-1}
\end{minipage}
\end{tabular}
\end{figure}

Regarding $\epsilon$ as a perturbation parameter, one writes $u(t)$ as a perturbation series
\begin{equation}\label{appl1:seri}
u(t)=u_0(t)+\epsilon u_1(t)+\epsilon^{2}u_2(t)+\epsilon^{3}u_3(t)+\cdots.
\end{equation}
Substituting Eq. (\ref{appl1:seri}) into Eq. (\ref{appl1:prob}), one obtains a perturbation solution
\begin{equation}\label{appl1:equ1}
u(t)=e^{-t}+\epsilon \left(e^{-t}-e^{-2t}\right)+\epsilon^{2}\left(\frac{1}{2}e^{-t}-2e^{-2t}+\frac{3}{2}e^{-3t}\right)+\cdots.
\end{equation}
The derivation of the perturbation solution (\ref{appl1:equ1}) at $t=0$ is
\begin{equation}\label{appl1:der-per}
u'(0)=-1+\epsilon-{\epsilon}^{2}+{\epsilon}^{3}-{\epsilon}^{4}+{\epsilon}^{5}-{\epsilon}^{6}+{\epsilon}^{7}-{\epsilon}^{8}+{\epsilon}^{9}-{\epsilon}^{10}+\cdots.
\end{equation}
Obviously, when $\epsilon\ge1$ the series (\ref{appl1:der-per}) is divergent, as shown in Fig. \ref{fig:appl1-1}.

Next, we analyse the problem by the traditional ADM. Following the procedure outlined in Section \ref{sect:TADM}, one obtains
\begin{equation}\label{appl1:equ2}
\begin{aligned}
&u_0(t)=1,\\
&u_1(t)=-t,\\
&u_2(t)=\epsilon t+\frac{1}{2}{t}^{2},\\
&u_3(t)=-{\epsilon}^{2}t-\frac{3}{2}\epsilon\,{t}^{2}-\frac{1}{6}{t}^{3},\\
\vdots
\end{aligned}
\end{equation}
The derivation of the ADM solution $u(t)$ at $t=0$ is
\begin{equation}\label{appl1:der-ADM}
u'(0)=-1+\epsilon-{\epsilon}^{2}+{\epsilon}^{3}-{\epsilon}^{4}+{\epsilon}^{5}-{\epsilon}^{6}+{\epsilon}^{7}-{\epsilon}^{8}+{\epsilon}^{9}-{\epsilon}^{10}+\cdots,
\end{equation}
which is the same as Eq. (\ref{appl1:der-per}) given by the perturbation method, as shown in Fig. \ref{fig:appl1-1}.

Finally, the ADMP is used to solve the problem (\ref{appl1:prob}). By means of (\ref{equ:partp}), one has
\begin{equation}\label{appl1:equ4}
\begin{aligned}
&v_0(t,c)=1,\\
&v_1(t,c)=-ct,\\
&v_2(t,c)=\left( c-1+{c}^{2}\epsilon \right) t+\frac{1}{2}{c}^{2}{t}^{2},\\
&v_3(t,c)=\left( 2c\epsilon-2{c}^{2}\epsilon-{c}^{3}{\epsilon}^{2} \right)t+ \left( c-{c}^{2}-\frac{3}{2}{c}^{3}\epsilon \right) {t}^{2}-\frac{1}{6}{c}^{3}{t}^{3},\\
&\vdots
\end{aligned}
\end{equation}
When $c=1$, they are the same as those in (\ref{appl1:equ2}) given by the ADM.

The $n$th-order approximation of $u'(0)$ is
\begin{equation}\label{appl1:der-appro}
\psi'_n(0,c)=\sum_{k=0}^{n}v'_k(0,c),
\end{equation}
where $v'_k(0,c),k=0,1,\cdots,n$ can be calculated via (\ref{appl1:equ4}) and the prime denotes differentiation with respect to $t$. For example, when $n=1,2,3,\cdots $, one has
\begin{equation*}
\begin{aligned}
&\psi'_1(0,c)=-c,\\
&\psi'_2(0,c)=-1+c^2\epsilon,\\
&\psi'_3(0,c)=-1+ \left( 2c-{c}^{2} \right) \epsilon-{c}^{3}{\epsilon}^{2},\\
&\vdots
\end{aligned}
\end{equation*}
It is easy to see that the $1st$-order ADMP approximation $\psi'_1(0,c)$ when $c=\frac{1}{1+\epsilon}$
coincides with the exact solution (\ref{appl1:der-exact}), whereas the $9$th-order perturbation
approximation and the $10$th-order ADM approximation diverge when $\epsilon\ge 1$, as shown in Fig. \ref{fig:appl1-1}.

\subsection{A nonlinear electrostatic cantilever NEMS model}

Beam-type electrostatic actuators have been widely used to construct nano-electromechanical system (NEMS). Over the last decades, the fundamental and applied researches as well as engineering developments in NEMS have undergone major improvements \cite{Duan2013, Ramezani2006}.

A beam-type NEMS actuator is modeled by a beam of length $L$, width $w$ with a uniform cross section and thickness $h$, which is suspended over a conductive substrate and separated by a dielectric spacer. In this paper, we discuss the cantilever NEMS. Based on the Euler-Bernoulli beam assumptions, the governing equation for a beam actuator with the first-order fringing field correction may be expressed as
\begin{equation}\label{appl2:equ1}
\frac{d^4Y}{dX^4}=\left(\frac{1}{(s-Y)^2}+\frac{0.65}{w(s-Y)}\right)F+\frac{1}{EI}F_k,\;k=3,4,
\end{equation}
where $Y, X, I, E$ and $F_k$ denote the deflection of the beam, the position along the beam as measured from the clamped end, the moment of inertia for the cross-sectional area of the beam, effective material modulus and intermolecular/quantum force per unit length of the beam respectively and $F=$ $\epsilon_0wV^2/2EI$ with $\epsilon_0=8.854\times10^{-12}C^{2} N^{-1} m^{-2}$ the dielectric constant of air, $s$ the original gap between the two electrodes without deflection and $V$ the applied voltage. The effective material modulus will become the plate modulus $E/(1-\nu^2)$, when the width satisfies $w\ge5h$.

The Van der Waals force \cite{Israelachvili2011} is
\begin{equation}\label{appl2:van}
F_3=\frac{Aw}{6\pi(s-Y)^3},
\end{equation}
with $A$ the Hamaker constant.
The Casimir force \cite{Lamoreaux2005} is
\begin{equation}\label{appl2:cas}
F_4=\frac{\pi^2\hbar vw}{240(s-Y)^4},
\end{equation}
with $\hbar=1.055\times10^{-34}J s$ the Planck's constant and $v$ the speed of light.
Substituting Eqs. (\ref{appl2:van}), (\ref{appl2:cas}) and the expressions $y=Y/s, x=X/L$ into Eq. (\ref{appl2:equ1}), yields
\begin{equation}\label{appl2:equ2}
\frac{d^4y(x)}{dx^4}=\frac{\alpha_k}{(1-y(x))^k}+\frac{\alpha_2}{(1-y(x))^2}+\frac{\alpha_1}{1-y(x)},\;k=3,4,
\end{equation}
where
\begin{equation*}
\alpha_4=\frac{\pi^2L^4hvw}{240EIs^5},\;\alpha_3=\frac{L^4Aw}{6\pi EIs^4}, \; \alpha_2=\frac{L^4\epsilon_0wV^2}{2EIs^3}, \;\alpha_1=\frac{0.65s}{w}\alpha_2.
\end{equation*}
The boundary conditions are
\begin{equation*}
y(0)=0, \;y'(0)=0, \; y''(1)=0, \; y'''(1)=0.
\end{equation*}
Setting $u(x)=1-y(x)$, then Eq. (\ref{appl2:equ2}) becomes
\begin{equation}\label{appl2:pro}
\frac{d^4u(x)}{dx^4}+\frac{\alpha_k}{u(x)^k}+\frac{\alpha_2}{u(x)^2}+\frac{\alpha_1}{u(x)}=0,
\end{equation}
with the boundary conditions
\begin{equation}\label{appl2:cons}
u(0)=1,\;u'(0)=0,\;u''(1)=0,\;u'''(1)=0.
\end{equation}

In this paper $k=3$ is considered. It is assumed that $\alpha_k=0.2, \alpha_2=0.5, \alpha_1=0.25$.
Following the procedure outline in Section \ref{sect:ADMP}, one has
\begin{equation}\label{appl2:equ3}
L[u(x)]+N[u(x)]=0,
\end{equation}
where
\begin{equation*}
L(\cdot)=\frac{d^4}{dx^4}(\cdot),\;N[u(x)]=\frac{\alpha_3}{u(x)^3}+\frac{\alpha_2}{u(x)^2}+\frac{\alpha_1}{u(x)}.
\end{equation*}
The inverse of the linear operator is taken to be
\begin{equation}\label{appl2:equ4}
L^{-1}=\int_0^{x}\int_0^{x_3}\int_1^{x_2}\int_1^{x_1}(\cdot)dtdx_1dx_2dx_3.
\end{equation}
The boundary conditions (\ref{appl2:cons}) then give
\begin{equation}\label{appl2:equ5}
\phi(x)=1.
\end{equation}

Combining Eqs. (\ref{appl2:equ3}), (\ref{appl2:equ4}) with (\ref{appl2:equ5}) and applying the recursion (\ref{equ:partp}), one obtains
\begin{equation*}
\begin{aligned}
v_0(x,c)=&1,\\
v_1(x,c)=&-c \left( {\frac {19{x}^{2}}{80}}-{\frac {19{x}^{3}}{120}}+{\frac{19{x}^{4}}{480}} \right),\\
v_2(x,c)=& \left( -{\frac {9139{c}^{2}}{288000}}-{\frac {19}{80}}+{\frac {19c}{80}} \right) {x}^{2}+ \left( {\frac {703\,{c}^{2}}{48000}}+{\frac {19}{120}}-{\frac {19c}{120}} \right) {x}^{3}\\
&+ \left( -{\frac {19}{480}}+{\frac {19c}{480}} \right) {x}^{4}-{\frac {703{c}^{2}{x}^{6}
}{576000}}+{\frac {703{c}^{2}{x}^{7}}{2016000}}-{\frac {703{c}^{2}{x}^{8}}{16128000}},\\
&\vdots
\end{aligned}
\end{equation*}
Then the $n$th-order approximation is
\begin{equation}\label{appl2:seri}
\psi_n(x,c)=\sum_{k=0}^{n}v_k(x,c).
\end{equation}
For example, when $n=1,2,\cdots$,
\begin{equation*}
\begin{aligned}
\psi_1(x,c)=&1-{\frac {19\,c{x}^{2}}{80}}+{\frac {19\,c{x}^{3}}{120}}-{\frac {19\,c{x}^{4}}{480}},\\
\psi_2(x,c)=&1- \left( {\frac {9139\,{c}^{2}}{288000}}+{\frac {19}{80}} \right) {x}^{2}+ \left( {\frac {703\,{c}^{2}}{48000}}+{\frac {19}{120}} \right){x}^{3}-{\frac {19\,{x}^{4}}{480}}-{\frac {703\,{c}^{2}{x}^{6}}{576000}}\\
&+{\frac {703\,{c}^{2}{x}^{7}}{2016000}}-{\frac {703\,{c}^{2}{x}^{8}}{16128000}},\\
&\vdots
\end{aligned}
\end{equation*}

The boundary value problem (\ref{appl2:pro})-(\ref{appl2:cons}) does not have a closed-form solution.
To determine whether the ADM solution (when $c=1$) is optimal or not, one uses the error remainder of
the $n$th-order approximation $\psi_n(x,c)$ defined by
\begin{equation}\label{appl2:err}
R_n(x,c)=L[\psi_n(x,c)]+N[\psi_n(x,c)],
\end{equation}
and the maximal error remainder of $\psi_n(x,c)$ defined by
\begin{equation}
MER_n=\max_{0\le x \le 1}|R_n(x,c)|.
\end{equation}

One takes $x_j = j/20 (j = 1, 2..., 20)$ as sample points and solves Eq. (\ref{equ:opti}). The optimal values of $c$
for $n=1$ to $n=10$ and the corresponding maximal error remainders are calculated, as shown in the second and third
columns of Table \ref{tab:appl2-1}, and the maximal error remainders of the ADM approximations (when $c=1$) are
shown in the fourth column. It is seen that the ADM approximations are not optimal and the ADMP approximations are
more accurate than the ADM approximations.
\begin{table}[ht]
\centering
\begin{tabular}{|c|cc|c|}\hline
$n$                &$c$                   &$MER_n(ADMP)$      &$MER_n(ADM)$(c=1)  \\
\hline
$1$                &1.10733               &1.01966E-1             &2.69754E-1  \\
$2$                &1.21693               &2.22708E-2             &1.01782E-1  \\
$3$                &1.21761               &2.99159E-3             &4.26177E-2  \\
$4$                &1.21455               &4.61594E-4             &1.89559E-2  \\
$5$                &1.21285               &7.58058E-5             &8.77839E-3  \\
$6$                &1.21203               &1.27843E-5             &4.18515E-3  \\
$7$                &1.21165               &2.17946E-6             &2.03977E-3  \\
$8$                &1.21147               &3.72997E-7             &1.01153E-3  \\
$9$                &1.21139               &6.38943E-8             &5.08700E-4  \\
$10$               &1.21135               &1.09425E-8             &2.58802E-4  \\
\hline
\end{tabular}
\caption{The $MER_n$ of the ADMP and the ADM }\label{tab:appl2-1}
\end{table}
\begin{figure}[ht]
\centering
\includegraphics[width=2in,height=2in]{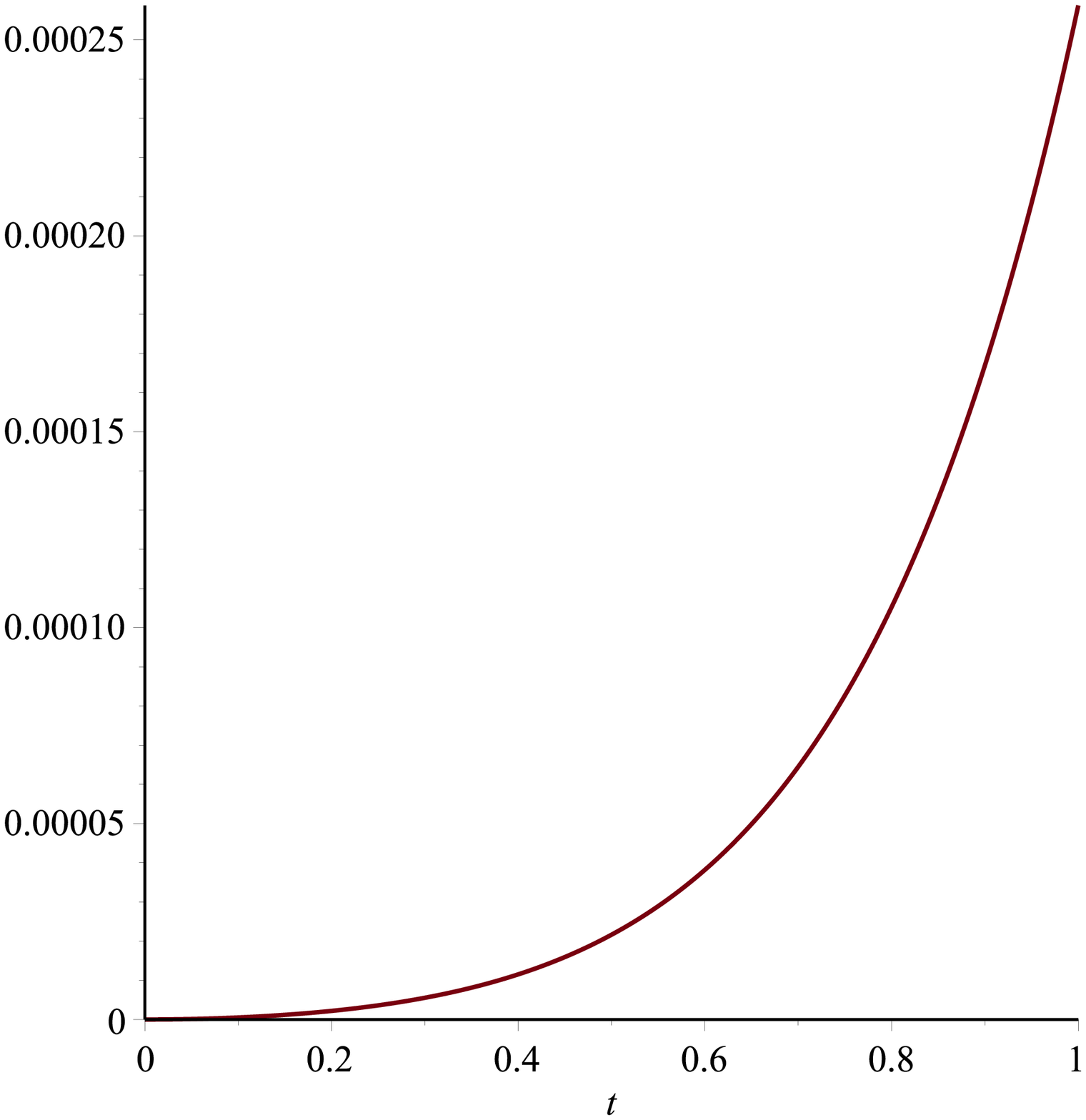}\includegraphics[width=2in,height=2in]{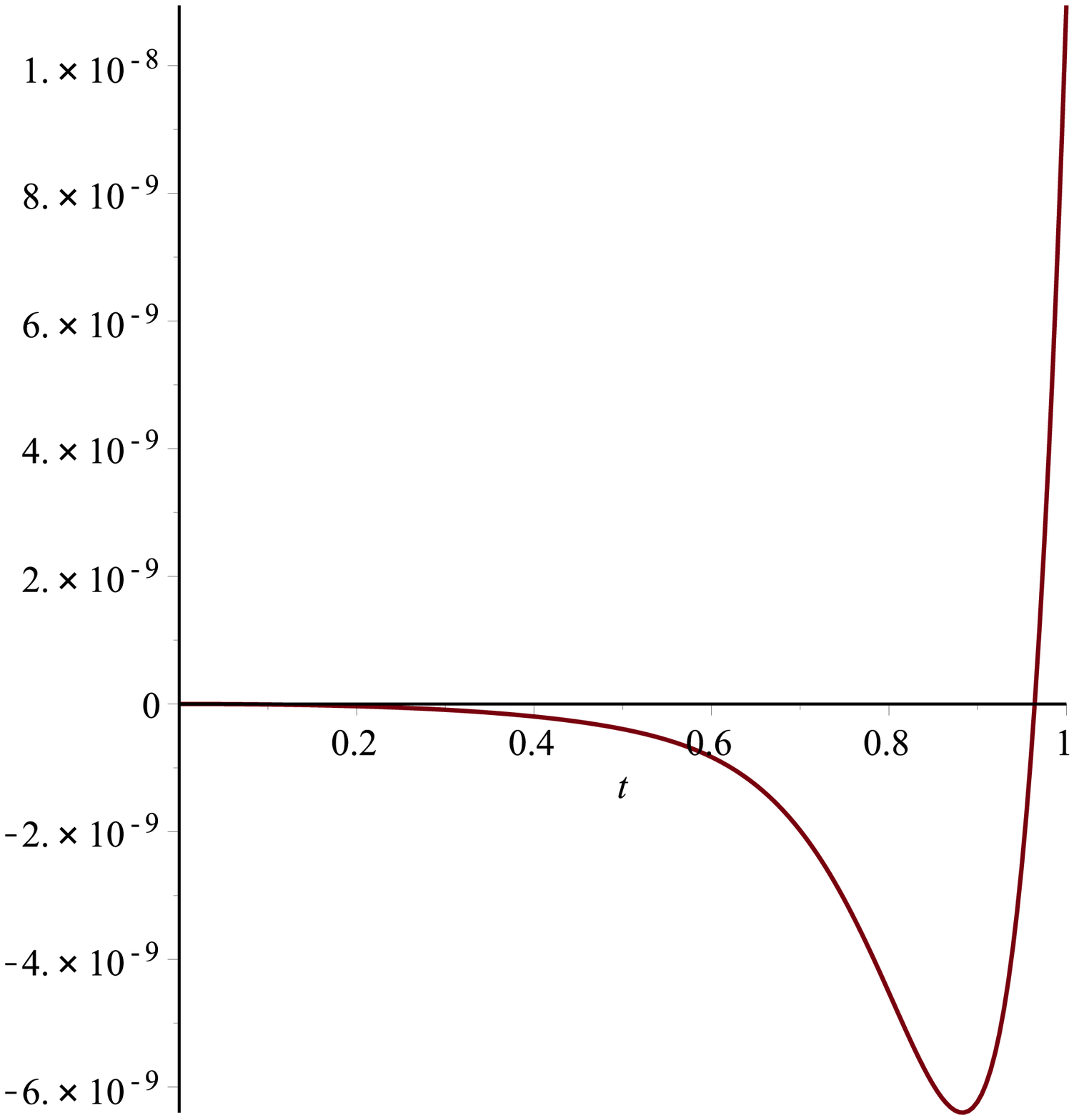}
\caption{LEFT: $R_{10}(t,1)$ for the ADM. RIGHT: $R_{10}(t,1.21135)$ for the ADMP.}\label{fig:appl2-1}
\end{figure}
To further demonstrate the accuracy of the ADMP approximation with the optimal value of $c$
over the ADM approximation, the error remainder $R_{10}(x,c)$
of the $10$th-order ADMP approximation when $c=1.21135$ and the error remainder $R_{10}(x,1)$ of the $10$th-order ADM approximation
are shown in Fig. \ref{fig:appl2-1}. It is seen that the error remainder of the ADMP approximation is much smaller than that of the
traditional ADM in the whole interval.

\subsection{A nonlinear Burgers' equation}\label{subsect:4.4}

Burgers' equation is one of the significant tools for describing nonlinear diffusion and dissipation phenomena. It is widely used to model the hydrodynamic motion \cite{Burgers1948}, soil-water flow \cite{Broadbridge1992} and nonlinear standing waves in constant-cross-sectioned resonators \cite{Bednarik2014}. In this part, a burgers' equation is considered with the form
\begin{equation}\label{appl3:problem}
\left\{
\begin{aligned}
&\frac{\partial u(x,t)}{\partial t}+u(x,t)\frac{\partial u(x,t)}{\partial x}+\frac{\partial^3u(x,t)}{\partial x^2\partial t}=0,\\
&u(x,0)=x.
\end{aligned}
\right.
\end{equation}
Eq. (\ref{appl3:problem}) has an exact solution
\begin{equation}\label{appl3:exact}
u(x,t)_{ex}=\frac{x}{1+t}.
\end{equation}

The linear operator is set as $L=\frac{d}{dt}(\cdot)$. Following the procedure outlined in Section \ref{sect:ADMP}, one has the series solution
\begin{equation}
v(x,t)=\sum_{k=0}^{+\infty}v_k(x,t,c),
\end{equation}
and an $n$th-order approximation
\begin{equation}\label{appl3:seri}
\psi_n(x,t,c)=\sum_{k=0}^{n}v_k(x,t,c),
\end{equation}
where $v_k(x,t,c),k=0,\cdots,n$ can be obtained via (\ref{equ:partp}). For example, when $n=1,2,\cdots,$
\begin{equation*}
\begin{aligned}
\psi_1(x,t,c)=& \left( 1-ct \right) x,\\
\psi_2(x,t,c)=&\left( 1-t+{c}^{2}{t}^{2} \right) x,\\
\psi_3(x,t,c)=&\left( 1-t-{c}^{2}{t}^{2}+2\,c{t}^{2}-{c}^{3}{t}^{3} \right) x,\\
\psi_4(x,t,c)=& \left( 1-t+{t}^{2}+2\,{c}^{3}{t}^{3}-3\,{c}^{2}{t}^{3}+{c}^{4}{t}^{4} \right) x,\\
\vdots
\end{aligned}
\end{equation*}

Consider the error of the $n$th-order approximation
\begin{equation}\label{appl3:error}
E_n(x,t,c)=\psi_n(x,t,c)-u(x,t)_{ex}.
\end{equation}

\begin{figure}[h]
\setcaptionwidth{2in}
\centering
\begin{tabular}{cc}
\begin{minipage}{2in}
\includegraphics[width=2in]{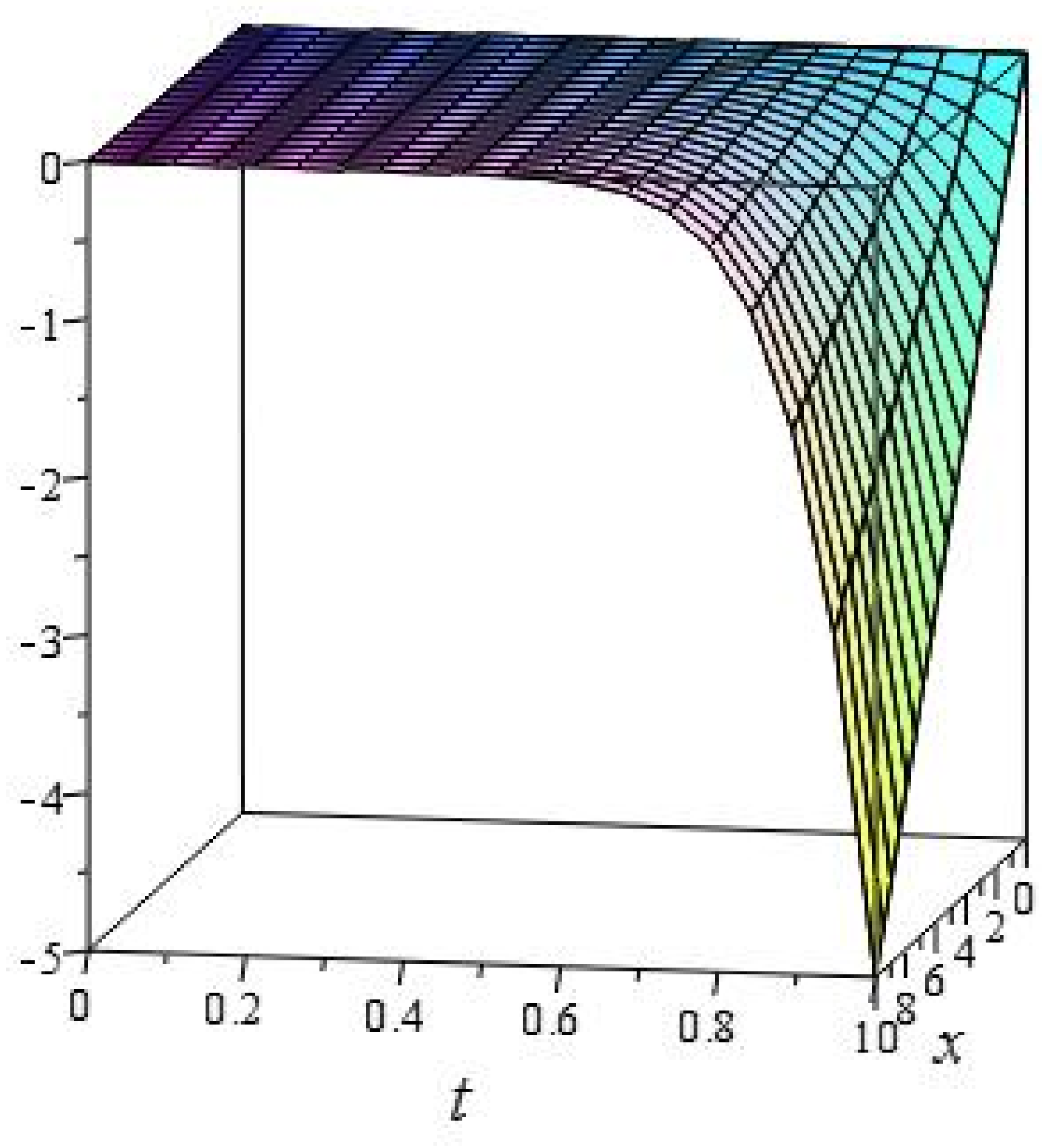}
\caption{$E_{10}(x,t,1)$ of the traditional ADM approximation for (\ref{appl3:problem}).}\label{fig:appl3-2}
\end{minipage}
\hspace{0.2in}
\begin{minipage}{2in}
\includegraphics[width=2in]{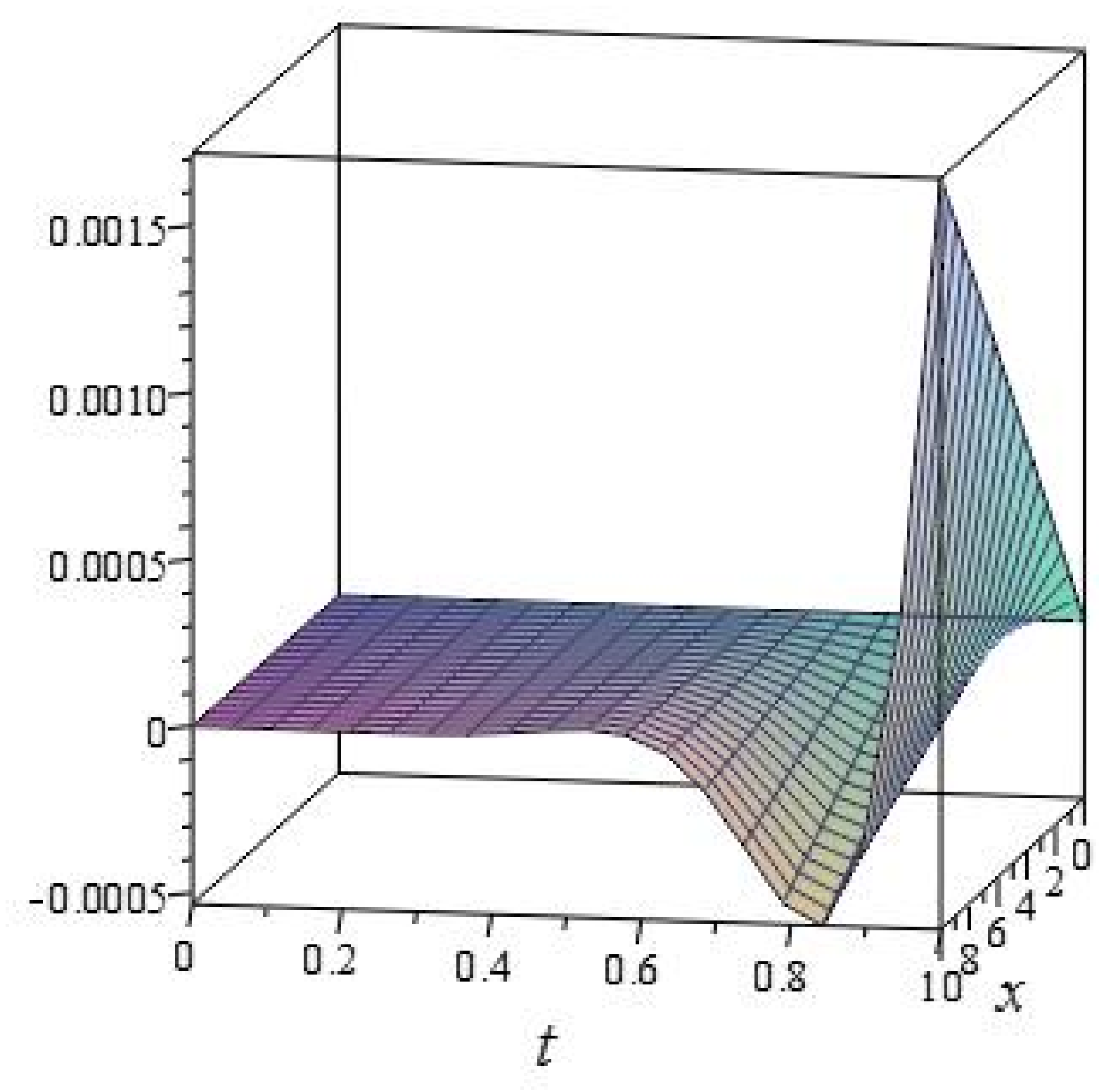}
\caption{$E_{10}(x,t,0.8)$ of the ADMP approximation for (\ref{appl3:problem}).}\label{fig:appl3-3}
\end{minipage}
\end{tabular}
\end{figure}

For $n=10$, Fig. \ref{fig:appl3-2} shows that the error $E_{10}(x,t,1)$ of the ADM approximation is very large. Therefore,
the ADM approximation is far away from the exact solution and it is essential to introduce a convergence-control parameter $c$
to make the results better.

By taking $20\times 20$ sample points $(x_i,t_j)=(i/2,j/20), i,j=1, 2, \ldots, 20$ and solving Eq. (\ref{equ:opti}),
one obtains the optimal value $c=0.8$. Fig. \ref{fig:appl3-3} shows that the error $E_{10}(x,t,0.8)$ of the ADMP approximation
is much smaller than the error $E_{10}(x,t,1)$ of the ADM approximation. To further demonstrate the accuracy of the ADMP approximation, the two approximation methods are contrasted in Fig. \ref{fig:appl3-1}. Therefore, the ADMP approximation is much more
accurate than the ADM approximation.

\subsection{A nonlinear RLW equation}
The RLW equation was first proposed by Peregrine for describing nonlinear dispersive waves. It has been used to model a large class of physics phenomena such as  the regularised long wave \cite{Peregrine1966}, nonlinear shallow wave and longitudinal dispersive waves in elastic rods. The RLW equation has attracted much attention. Existence and uniqueness of the solution to the RLW equation has been given in \cite{Bona1973}. Several numerical method such as finite-difference method, a lumped Galerkin method and septic splines have been introduced in the literature \cite{Bhardwaj2000}. Consider the RLW equation with this form
\begin{equation}\label{appl4:prob}
\left\{
\begin{aligned}
&\frac{\partial u(x,t)}{\partial t}+\frac{\partial u(x,t)}{\partial x}+u(x,t)\frac{\partial u(x,t)}{\partial x}+\frac{\partial^3u(x,t)}{\partial x^2\partial t}=0,\\
&u(x,0)=x.\\
\end{aligned}
\right.
\end{equation}
The problem (\ref{appl4:prob}) has an exact solution
\begin{equation}\label{appl4:exact}
u(x,t)_{ex}=\frac{x-t}{1+t}.
\end{equation}

\begin{figure}[h]
\setcaptionwidth{2in}
\centering
\begin{tabular}{cc}
\begin{minipage}{2in}
\includegraphics[width=2in]{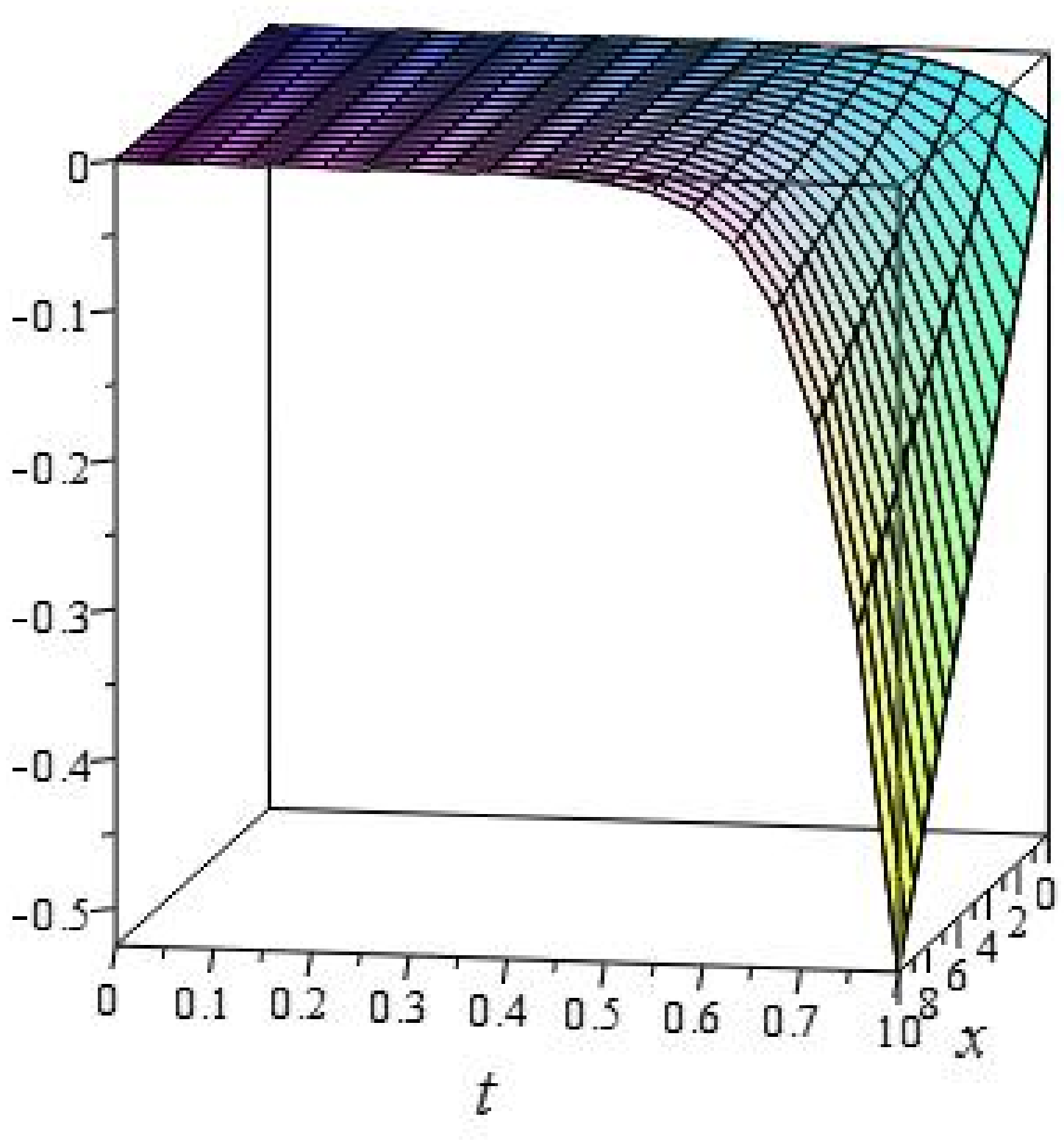}
\caption{$E_{10}(x,t,1)$ of the traditional ADM approximation for (\ref{appl4:prob}).}\label{fig:appl4-1}
\end{minipage}
\hspace{0.2in}
\begin{minipage}{2in}
\includegraphics[width=2in]{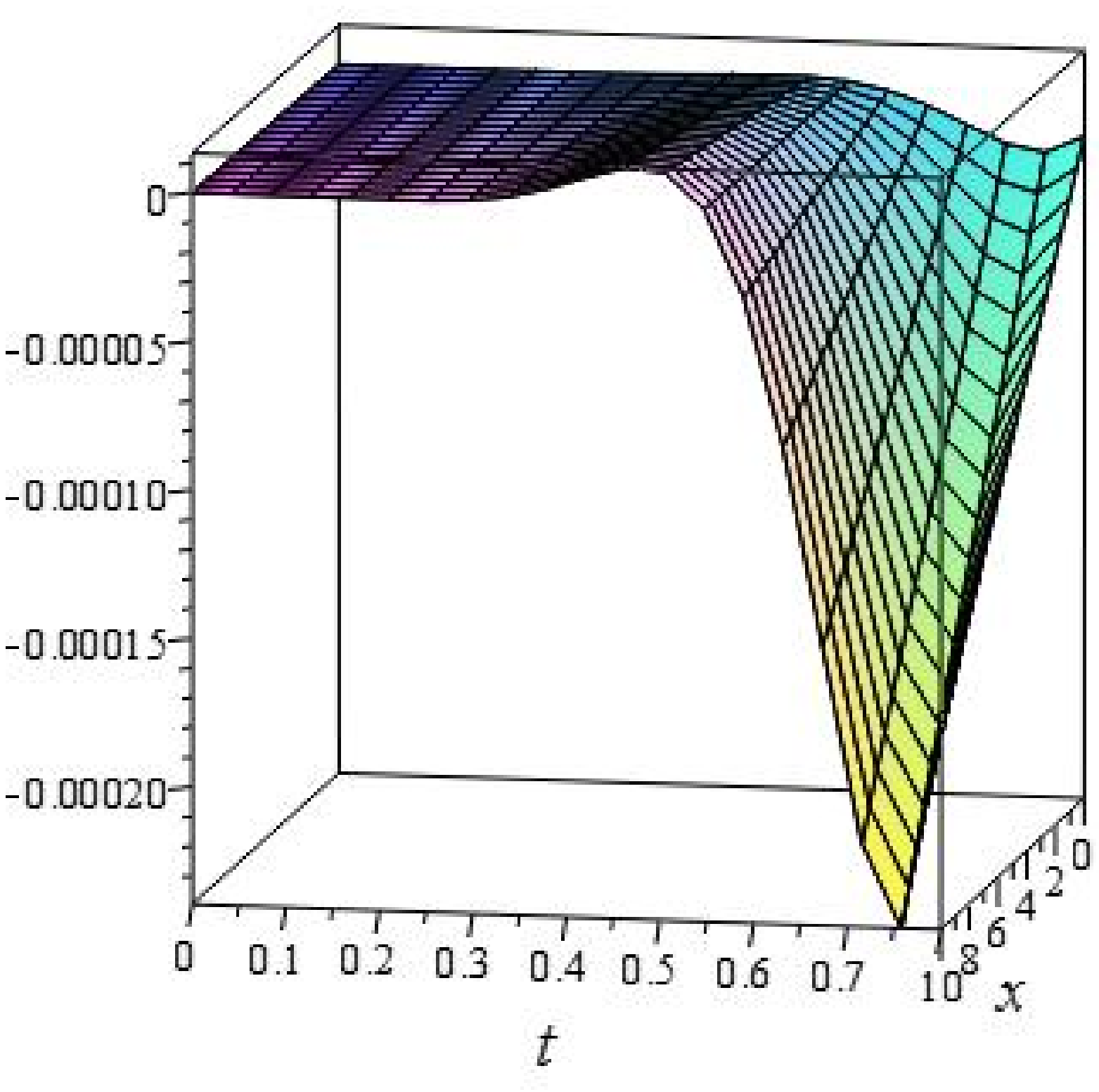}
\caption{$E_{10}(x,t,0.81)$ of the ADMP approximation for (\ref{appl4:prob}).}\label{fig:appl4-2}
\end{minipage}
\end{tabular}
\end{figure}

Following the procedure outlined in Section \ref{sect:ADMP}, one obtains the series solution
\begin{equation}
v(x,t)=\sum_{k=0}^{+\infty}v_k(x,t,c),
\end{equation}
and an $n$th-order approximation
\begin{equation}
\psi_n(x,t,c)=\sum_{k=0}^{n}v_k(x,t,c)
\end{equation}
to the problem (\ref{appl4:prob}). For example, when $n=1,2,\cdots,$ one has
\begin{equation*}
\begin{aligned}
\psi_1(x,t,c)&=x-ct-cxt,\\
\psi_2(x,t,c)&=x-t-xt+{c}^{2}x{t}^{2}+{c}^{2}{t}^{2},\\
\psi_3(x,t,c)&=x-t-xt-{c}^{2}x{t}^{2}-{c}^{2}{t}^{2}+2\,cx{t}^{2}+2\,c{t}^{2}-{c}^{3}x{t}^{3}-{c}^{3}{t}^{3},\\
&\vdots
\end{aligned}
\end{equation*}

Similar to the problem in Section \ref{subsect:4.4}, one obtains the optimal value $c=0.81$ for the $10$th-order
ADMP approximation. Fig. \ref{fig:appl4-1} and Fig. \ref{fig:appl4-2} show that the error $E_{10}(x,t,0.81)$ of the ADMP approximation
is much smaller that the error $E_{10}(x,t,1)$ of the ADM approximation. Therefore, the ADMP approximation is much more
accurate than the ADM approximation.

\section{Conclusion}
We have introduced a new convergence-control parameter into the traditional ADM, which can be used to adjust and control the convergence region and rate of the resulting series solution. Furthermore, we have proposed an approach for finding the optimal value of the convergence-control parameter. The traditional ADM is only a special case of the ADMP when $c=1$. Usually, the $n$th-order approximation $\psi_n(x,t,c)$ is not optimal when $c=1$. The ADMP is significantly helpful for obtaining more accurate approximation.

\section*{Acknowledgement}

The authors would like to thank Professor Junsheng Duan (College of Science, Shanghai Institute of Technology) for his suggestions.


\begin{thebibliography}{1}

\bibitem{Nayfeh2000} A. H. Nayfeh, \textit{Perturbation Methods}, John Wiley\&Sons, New York, 2000.

\bibitem{Lyapunov1992} A. M. Lyapunov, \textit{General Problem on Stability of Motion}, Taylor\&Francis, London, 1992.

\bibitem{Karmishin1990} A. V. Karmishin, A. T. Zhukov and V. G. Kolosov, \textit{Methods of Dynamics Calculation and Testing for Thin-Walled Structures}, Moscow, Mashinostoyenie, 1990.

\bibitem{Adomian1983} G. Adomian and R. Rach, Nonlinear stochastic differential-delay equations, J. Math. Anal. Appl.
\textbf{91}(1) (1983), 94--101.

\bibitem{Deeba1995} E. Y. Deeba and S. A. Khuri, A Decomposition method for solving the nonlinear Klein-Gordon equation, J. Comput. Phys. \textbf{124}(0071) (1995), 442--448.

\bibitem{Wazwaz2014} A. M. Wazwaz, R. Rach and J. S. Duan, A study on the systems of the Volterra integral forms of the Lane¨CEmden equations by the Adomian decomposition method, Math. Meth. Appl. Sci. \textbf{37}(1) (2014):10-19.

\bibitem{Ebaid2008} A. Ebaid, A new numerical solution for the MHD peristaltic flow of a bio-fluid with variable viscosity
   in a circular cylindrical tube via Adomian decomposition method, Phys. Lett. A. \textbf{372}(32) (2008), 5321--5328.

\bibitem{Singla2015} R. K. Singla and R. Das, Adomian decomposition method for a stepped fin with all temperature-dependent modes of heat transfer, Int. J. Heat Mass Tran. \textbf{82} (2015), 447--456.

\bibitem{Siddiqui2010} A. M. Siddiqui, M. Hameed, B. M. Siddiqui and Q. K. Ghori, Use of Adomian decomposition method in the study of parallel plate flow of a third grade fluid, Commun. Nonlinear Sci. \textbf{15}(9) (2010), 2388--2399.

\bibitem{Adomian1992} G. Adomian and R. Rach, Noise terms in decomposition series solution, Comput. Math. Appl. \textbf{24}(11) (1992), 61--64.

\bibitem{Rach1987} R. Rach, On the Adomian (decomposition) method and comparisons with Picard's method, J. Math. Anal. Appl. \textbf{128}(2) (1987), 5321--5328.

\bibitem{Kaya2004} D. Kaya, A reliable method for the numerical solution of the kinetics problems, Appl. Math. Comput. \textbf{156}(1) (2004), 261--270.

\bibitem{Wazwaz1999} A. M. Wazwaz, A reliable modification of Adomian decomposition method, Appl. Math. Comput. \textbf{102}(1) (1999), 261--270.

\bibitem{El-Kalla2011}I. L. El-Kalla, Error estimate of the series solution to a class of nonlinear fractional differential equations, Commun. Nonlinear Sci. \textbf{16}(3) (2011), 1408-1413.

\bibitem{Duan2011} J. S. Duan and R. Rach, A new modification of the Adomian decomposition method for solving boundary value problems for higher order nonlinear differential equations, Appl. Math. Comput. \textbf{218}(8) (2011), 4090-4118.

\bibitem{Lin2014} Y. W. Lin, T. T. Lu and C. K. Chen, Large interval solution of the Emden-Fower equation using a modified Adomian Decomposition Method with an integrating factor, Anziam J. \textbf{56}(2) (2014), 192--208.

\bibitem{Zhang2015} X. L. Zhang and S. X. Liang, Adomian decomposition method is a special case of Lyapunov's artificial small parameter method, Phys. Lett. A \textbf{48} (2015), 177--179.

\bibitem{Duan2010} J. S. Duan, Recurrence triangle for Adomian polynomials, Appl. Math. Comput. \textbf{216}(4) (2010), 1235-1241.

\bibitem{Lu2014} L. Lu and J. S. Duan, How to select the value of the convergence parameter in the Adomian decomposition method, CMES-Comp. Model. Eng. Sci. \textbf{97}(1) (2014), 35-52.

\bibitem{Duan2013} J. S. Duan, R. Rach and A. M. Wazwaz, Solutin of the model of beam-type micor- and nano-scale electrostatic actuators by a new modified Adomian decomposition method for nonlinear boundary value problems, Int. J. Nonlin. Mech. \textbf{49} (2013), 159-169.

\bibitem{Liang2009} S. X. Liang and D. J. Jeffrey, Comparison of homotopy analysis method and homotopy perturbation method through an evolution equation, Commun. Nonlinear Sci. \textbf{14} (2009), 4057--4064.

\bibitem{Ramezani2006} A. Ramezani, A. Alasty and J. Akbari, Closed-form solutions of the pull-in instability in nano-cantilevers under electrostatic and intermolecular surface forces, Int. J. Solids Struct. \textbf{44}(14--15) (2006), 4925--4941.

\bibitem{Israelachvili2011} J. N. Israelachvili, \textit{Intermolecular and Surface Forces}, Academic press, London, 2011.

\bibitem{Lamoreaux2005} S. K. Lamoreaux, The Casimir force: background, experiments, and applications, Rep. Prog. Phys. \textbf{68}(1) (2005), 201--236.

\bibitem{Burgers1948} J. M. Burgers, A mathematical model illustrating the theory of turbulence, Adv. Appl. Mech. \textbf{1} (1999), 171--199.

\bibitem{Broadbridge1992} P. Broadbridge, R. Srivastava and T.-C. Jim Yeh, Burgers' equation and layered media: exact solutions and applications to soil-water flow, Math. Comput. Model. \textbf{16}(11) (1992), 163--169.

\bibitem{Bednarik2014} M. Bednarik and M. Cervenka, Equations for description of nonlinear standing waves in constant-cross-sectioned resonators, J. Acoust. Soc. Am. \textbf{135}(3) (2014), 135--139.

\bibitem{Peregrine1966} D. H. Peregrine, Calculations of the development of an undular bore, J. Fluid Mech. \textbf{25}(2) (1966), 321--330.

\bibitem{Bona1973} J.L. Bona and P.J. Bryant, A mathematical model for long waves generated by wave makers in nonlinear dispersive systems, Proc. Cambridge Philos. Soc. \textbf{73}(2) (1973), 391--405.

\bibitem{Bhardwaj2000} D. Bhardwaj and R. Shankar, A computational method for regularized long wave equation,Comput. Math. Appl. \textbf{40}(12) (2000), 1397--1404.










\end{thebibliography}
\end{document}